\documentclass[a4paper, english, 11pt]{amsart}
\usepackage[top = 3cm, bottom = 3cm, left = 3cm, right = 3cm]{geometry}
\usepackage[english]{babel}
\usepackage{lmodern}

\usepackage{amssymb}
\usepackage{amsmath}
\usepackage{hyperref}

\usepackage{listings}
\usepackage{xcolor}

\usepackage{import}
\usepackage{xifthen}
\usepackage{pdfpages}
\usepackage{transparent}

\newtheorem*{thma}{Theorem A}
\newtheorem*{thmb}{Theorem B}
\newtheorem{thm}{Theorem}[section]

\newtheorem{lemma}[thm]{Lemma}
\newtheorem{corollary}[thm]{Corollary}

\theoremstyle{definition}
\newtheorem{definition}{Definition}[section]
\newtheorem{remark}{Remark}[section]

\newtheorem{obs}{Remark}

\newcommand{\A}{\mathbb{A}}
\newcommand{\R}{\mathbb{R}}

\title{Rotational Chaos in the driven Pendulum} 

\author{Maciej J. Capi\'nski}
\address{AGH University of Krak\'ow, Krak\'ow, Poland; maciej.capinski@agh.edu.pl}
\author{Sof\'ia Llavayol}
\address{Facultad de Ciencias Económicas y de Administración, UdelaR, Montevideo Uruguay; sofia.llavayol@fcea.edu.uy.}
\author{Alejandro Passeggi}
\address{Facultad de Cienicas, UdelaR, Montevideo Uruguay; apasseggi@cmat.edu.uy.}

\begin{document}

\begin{abstract}
Building on recent developments in the topological theory of surface dynamics, we develop readily applicable criteria for computer-assisted proofs of the existence of rotational horseshoes and apply them to the \emph{Driven Pendulum}.
These tools also allow us to localize the resulting \emph{instability region}, and to compute a lower bound for the length of their rotation sets. This, in turn, implies the existence of periodic points realizing a wide range of rational rotation numbers. The obtained results correspond to physically meaningful parameter values of the Driven Pendulum. 
The method also applies to periodically forced Hamiltonian systems whenever the stated conditions are satisfied.

\end{abstract}

\maketitle


\section{Introduction}

Although chaotic behavior has been extensively studied in specific systems arising from physics, engineering, and other scientific disciplines, its rigorous mathematical validation often involves substantial technical challenges.
Moreover, when such results are available, their proofs are often 
tailored to the particular system under consideration and apply only within restricted parameter regimes, sometimes far removed from physically relevant ones. Frequently, they rely on perturbations of integrable or essentially one-dimensional models. At the same time, much of the mathematical theory of dynamical systems has focused on developing and understanding abstract models of complex behavior---such as Hyperbolic systems and their generalizations---rather than on proving that these structures occur in concrete systems arising in other scientific contexts.
This can be regarded as a serious and persistent issue
in the theory, particularly in light of Poincar\'e's
reflections on the role of mathematics in other sciences
when explicit integration of differential equations is
no longer attainable \cite[pp.~36--37]{poincare_method},
and Smale's remarks surrounding his celebrated
Problem 14 \cite[p.~17]{smale_problems}.


In non-rigorous treatments, for instance, chaos is commonly diagnosed through numerical indicators, most notably Lyapunov exponents. Such quantities necessarily capture only part of the underlying dynamics. Since Smale's foundational work clarified the mathematical nature of chaotic behavior, it has been understood that chaos in low-dimensional regimes is more appropriately associated with a rich, finitely generated dynamical structure, now commonly referred to as a \emph{horseshoe} \cite{SmaleRio}. A horseshoe not only yields numerical manifestations of chaos, such as positive entropy or positive Lyapunov exponents under suitable regularity assumptions, but also provides detailed information about the organization of the orbits.

A given system may contain many distinct horseshoes, as by the well-developed theory of generic area-preserving dynamics, see the recent developments \cite{LeCalvezSambarinoHomoclinic,LeCalvezSambarinoNonContractible} and reference therein.
However, as emphasized in \cite{hubbard_pendulum}, when studying a specific system, the horseshoes of interest are not arbitrary: one seeks horseshoes encoding dynamically meaningful information related to the topology and geometry of the particular configuration space.
For example, when the configuration space is a circle, it is natural to look for horseshoes carrying rotational information around the circle, rather than horseshoes whose orbits remain confined near an equilibrium point and never wind around the configuration space. Such structures are known in the topological theory of surface dynamics as \emph{rotational horseshoes}. At the end of the introduction we discuss the physical meaning of this
dynamical structure in the particular example of the \emph{driven pendulum}.

Recent advances in the topological theory of surface dynamics have led to simple criteria guaranteeing the existence of rotational horseshoes \cite{passeggi_tal_conditions,passpir}. Roughly, one identifies open sets $U_0,U_1$ in the annulus with a prescribed \emph{rotational difference}; if they visit one another under the dynamics, a rotational horseshoe follows. These criteria have already found applications in the theoretical setting and yielded strong conclusions for well-known families of concrete systems \cite{passeggi_tal_groeger_capinski}.

\smallskip

The motivation for this work is to demonstrate the applicability of these topological criteria to systems arising in physics. We first derive a criterion for periodically forced Hamiltonian systems (Theorem 3.4), and then apply it to the driven pendulum, a classical model closely connected with the development of dynamical systems, from Poincaré's pioneering investigations of homoclinic phenomena \cite{poincare1890} to Franks' celebrated generalization of the Poincaré–Birkhoff theorem \cite{franks}.

We obtain a simple computer-assisted proof of the existence of rotational horseshoes for physically relevant parameter values, going beyond the near-integrable regimes typically required by classical perturbative approaches, such as the Melnikov method. Moreover, our results provide rigorous quantitative information: the localization of an instability region and a lower bound for the width of its rotation set.

Beyond the \emph{driven pendulum}, our criterion given by Theorem \ref{t.diffeq} can be applied to return maps of
periodically forced Hamiltonian systems, widely studied
in many models arising in physics and engineering, see for instance,
\cite{chirikov1979,rechester1980,mackay1984,balescu1998,
balescu1998nontwist,abdullaev2004,hensinger2001,
hensinger2001b}, and we hope that the present work
will encourage further applications in these directions.

\smallskip

Let us introduce a first instance of this application. Recall that the driven-pendulum is governed by
\begin{equation}
    \label{eq-driven-pendulum}
    \ddot{x}
    =
    -\frac{g}{\ell}\sin(x)+A\sin(\omega_0t).
\end{equation}
Introducing the auxiliary variable $y=\dot{x}$, we obtain a non-autonomous ordinary differential equation, periodic in time with period $2\pi/\omega_0$, on the annulus
\begin{equation}
    \label{eq:annulus}
    \A=\R^2/{\sim},
    \qquad
    (x_1,y_1)\sim(x_2,y_2)
    \quad\Longleftrightarrow\quad
    x_2-x_1\in2\pi\mathbb Z \mbox{ and }y_1=y_2.
\end{equation}
The corresponding stroboscopic map, obtained by following the flow for one forcing period $2\pi/\omega_0$,
\[
f:\A\longrightarrow\A,
\]
is a homeomorphism homotopic to the identity. 

Since the differential equation has zero divergence in the variables $(x,y)$, the map $f$ preserves the standard Lebesgue measure on $\A$, naturally induced from the quotient. It can also be regarded as a Poincar\'e return map in the extended phase space. Our first result is not intended as an exhaustive exploration of the parameter space. Instead, it illustrates the amount of rigorous dynamical information that can be extracted from a single representative example. The parameters were chosen within a reasonable physical range and so that the orbit segments connecting the relevant sets $U_0$ and $U_1$ are as short as possible.

\begin{thma}
    For the parameter values
    \[
        g=1.5,
        \qquad
        A=0.5,
        \qquad
        \ell=1,
        \qquad
        \omega_0=1,
    \]
    the associated Poincar\'e map $f$ admits a rotational horseshoe.

    Moreover, this horseshoe is contained in an instability region
    $\mathcal I$ that intersects the projections of the boxes
       \begin{equation}
        \label{eq:th-boxes}
        [2.8916,3.2231]\times[0.4831,0.5953]
        \quad\text{and}\quad
        [1.7910,2.0412]\times[-2.6007,-2.3302].
    \end{equation}
    The length of the rotation interval
    $\rho_{\mathcal I}(F)$ is at least $2$.

    In particular, one can choose a lift $F$ of $f$ such that
    \[
        \rho_{\mathcal I}(F)\supset[0,2].
    \]
    For every rational number $p/q\in[0,2]$, there exists a periodic
    orbit of $f$, contained in $\mathcal I$, with rotation number
    $p/q$.
\end{thma}

Let us list some remarks about the result.

\begin{obs}
In accordance with the preceding discussion of the criteria used, the resulting computer-assisted proof is simple and can be readily adapted to other parameter values, making a systematic exploration of the parameter space feasible, see Theorem B for a first instance of such exploration.
Let us emphasize again that we want to use this example in order to show that it can be easily applied in different instances where these kind of mapsappears, as for instance those mentioned in the introduction.
\end{obs}

\begin{figure}
    \centering
    \includegraphics[width=\textwidth]{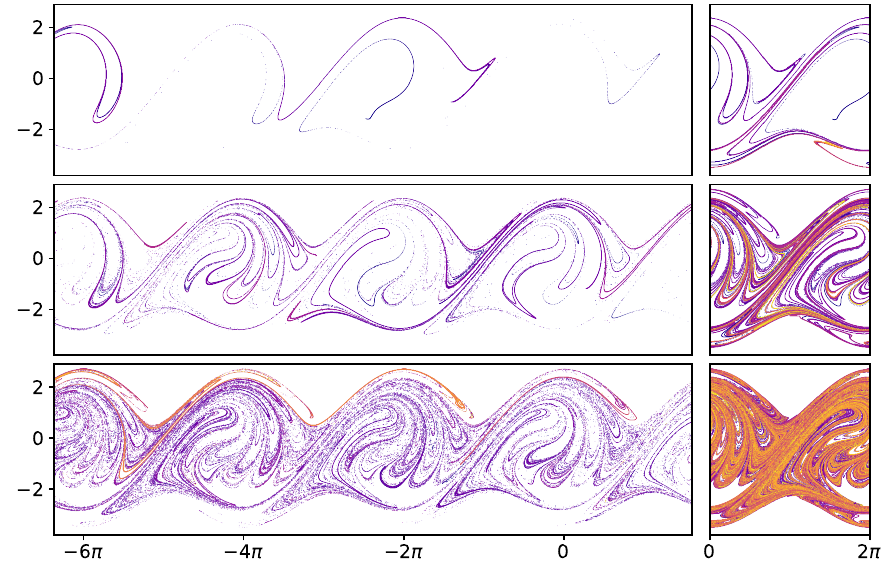}
    \caption{
        Numerical simulations of $F^n(U_1)$ for
        $n=3$, $6$, and $10$, respectively, where $F$ is the lift
        considered in Theorem~A and $U_1$ is the second box defined in \eqref{eq:th-boxes}.
        The final image provides a numerical approximation of the
        instability region containing the rotational horseshoe
        established in Theorem~A as mention in Remark \ref{r.instreg}.
    }
    \label{chaotic-sea}
\end{figure}

\begin{obs}\label{r.instreg}

The classical notions of an \emph{instability region} and of the \emph{rotation set} $\rho_K(F)$ of an $f$-invariant set $K\subset\A$, computed with respect to a lift $F$ of $f$, are recalled in Section~\ref{s.preliminaries}. In view of the information provided by Theorem A, we propose the following heuristic procedure to approximate the instability region. The theoretical picture suggests looking for pairs of points, one in $U_0$ and the other in $U_1$, whose forward orbits approach each other. Since such asymptotic behavior cannot be verified numerically, we replace it by a finite-time criterion: for a sufficiently large fixed iterate $N$, we select pairs whose images at time $N$ are $\varepsilon$-close. This informally means that such a pair can be regarded as lying in the same stable manifold of a point in the instability region, and its orbit therefore provides information about that region. We then iterate the selected points in $U_1$ (or those in $U_0$) to obtain the desired approximation. For simplicity, Figure~\ref{chaotic-sea} shows a simulation of the iterates of all points in $U_1$.


\end{obs}

\begin{obs}

The last part of the statement shows the existence of infinitely many periodic solution for our system, showing different rotational behavior. This
of course, has a direct physical impact, as the meaning of the rotation number of a periodic orbit is transparent: It measures the average angular displacement of the orbit in the universal cover per forcing period. After the corresponding normalization by the forcing period, it therefore determines the mean angular velocity. In particular, a nonzero rotation number implies that the trajectory winds around the circular configuration space instead of remaining confined near an equilibrium.

\end{obs}

As mention before, we also demonstrate that our method is applicable for automated parameter exploration, combined with a computer assisted proof. We have validated the following result.

\begin{thmb}
For the parameter values
    \[
        g=1.0,
        \qquad
        A\in [1.0,1.5],
        \qquad
        \ell=1,
        \qquad
        \omega_0=1,
    \]
    for over $96.55\%$ of the parameters $A$ in the interval $[1.0,1.5]$ the associated Poincar\'e map $f$ admits a rotational horseshoe.
\end{thmb}

\begin{obs}
We separate the two results in Theorems A and B because the proof of Theorem B uses the fact that the system admits periodic orbits and uses symmetric sets. In contrast, the proof of Theorem A does not rely on any such properties.
\end{obs}

\subsection{Physical interpretation of Rotational Horeseshoes}

In this paragraph we briefly explain the physical significance of a rotational horseshoe. This discussion follows the perspective developed in \cite{hubbard_pendulum}.

\begin{figure}
    \centering
    \def\svgwidth{\columnwidth}
    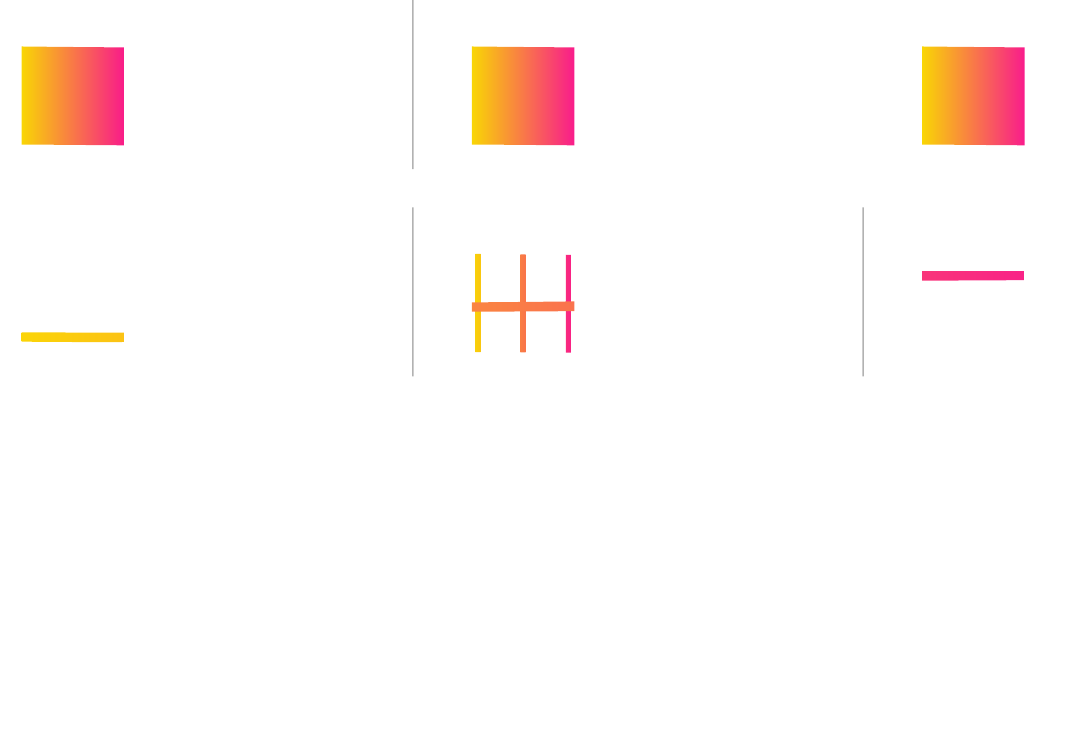

    \caption{
        A three-legged rotational horseshoe. The upper figure shows
        a rectangle $\widetilde R$ and its image $F(\widetilde R)$,
        which intersect in the required topological configuration,
        together with the consecutive translates
        $\widetilde R_{-1}$ and $\widetilde R_1$. Here
        \[
            \widetilde R_j
            :=
            \widetilde R+(2\pi j,0),
            \qquad j\in\mathbb Z.
        \]
        The lower figure shows the corresponding configuration after
        projection onto the annulus $\A$.
    }
    \label{f.horseshoe}
\end{figure}

The defitnition of rotational horseshoes is inspired by Smale's original horseshoe, but is adapted to a setting in which global hyperbolicity cannot be expected. Accordingly, we work with a $C^0$ version of a horseshoe in which the chaotic behavior is reflected by orbits that wind around the annular phase space. Figure~\ref{f.horseshoe} illustrates the type of rotational horseshoe considered here. We refer to Section~\ref{s.preliminaries} for the precise definition and to
\cite{GuiheneufMiliton2025,PasseggiPotrieSambarino2018,LeCalvezTal2022}
for related results.

Such an object allows us to control the behavior of solutions without requiring explicit formulas, a paradigm that has largely disappeared from the study of differential equations since the emergence of chaos theory. A rotational horseshoe admits a Markov-type symbolic coding. For a $k$-legged horseshoe, the dynamics can be encoded using $k$ symbols, with every admissible sequence corresponding to at least one orbit---and to exactly one orbit in the uniformly hyperbolic case. Crucially, this symbolic coding also records how the corresponding physical trajectories wind around the annulus.

For example, suppose that the driven pendulum has forcing period $T_0$ and that its stroboscopic map admits a three-legged rotational horseshoe of the type represented in Figure~\ref{f.horseshoe}, encoded by the symbols $-1$, $0$, and $1$. Consider the finite word
\begin{center}
    \large\textbf{%
        $1\quad -1\quad 0\quad 0\quad 1\quad 1\quad -1$}.
\end{center}
The symbolic dynamics guarantees the existence of an orbit whose behavior during the successive time intervals between multiples of $T_0$ consists of
\begin{center}
    \textbf{one counterclockwise turn},
    followed by
    \textbf{one clockwise turn},
    then two
    \textbf{homotopically trivial displacements},
    followed by
    \textbf{two counterclockwise turns},
    and finally
    \textbf{one clockwise turn}
\end{center}
around the circular configuration space.

Each of the expressions highlighted above admits a precise mathematical interpretation through the definition of a rotational horseshoe given in Section~\ref{s.preliminaries}. Moreover, since the finite word can be repeated indefinitely, the topological theory of surface dynamics guarantees the existence of a periodic solution of the differential equation exhibiting precisely this prescribed rotational behavior
\cite{GuiheneufMiliton2025}.

This example illustrates why detecting \emph{rotational horseshoes} in physical systems is particularly significant. Such a horseshoe does not merely establish the existence of chaos; it describes how chaotic trajectories interact with the global topology of the configuration space. In this sense, it carries substantially more system-specific information than the mere detection of a homoclinic intersection associated with a periodic orbit. We refer again to the illuminating discussion of J.~Hubbard in \cite{hubbard_pendulum}. We believe that identifying such structures can provide valuable insight into the systems under consideration. This perspective may be viewed as a response, in the spirit of the ideas expressed on page~4 of Poincar\'e's essay \cite{poincare_avenir}, to the problem of extracting qualitative information about solutions even when explicit formulas are unavailable.

\smallskip

The rigorous detection of Smale horseshoes has a long and rich history. One classical approach is to establish the existence of transverse homoclinic or heteroclinic intersections between the stable and unstable manifolds of suitable invariant objects and then apply the Birkhoff--Smale theorem. Such intersections can be obtained either through perturbative arguments or by computer-assisted proofs. The latter make it possible to validate transverse intersections at physically relevant parameter values; see, for example,
\cite{Cap2012,Jay2023,Jay2019-2,Jay2019,Jay2014}.

A second approach is to verify directly the topological alignment of the sets involved in a horseshoe construction. This method is also well suited to computer-assisted proofs; see, for instance,
\cite{Galias1998,Galias2001,Wilczak2003-2,Wilczak2006,Wilczak2003,Wilczak2009,Wilczak2020}.
Both approaches typically require a detailed analysis of the specific system: either an explicit construction and validation of intersections between invariant manifolds or the identification of suitable domains forming a horseshoe.

As mentioned above, recent developments provide a different route. The criteria employed in the present work are purely topological and require no control of derivatives, invariant manifolds, or previously located periodic orbits. This makes their computer-assisted implementation considerably more direct and opens the possibility of applying the same methodology systematically across broad parameter families. 

The method can be applied without a detailed analysis of the system's behaviour. The sets it requires can be found using simple `blind shooting' algorithms, which worked well in the example we investigated. Because the method does not rely on computing invariant manifolds and their intersections, it yields proofs that are fast in terms of computation time. It also does not require long integration times, except potentially between two distinct points belonging to two distinct sets. This is an advantage both for the speed of the method and for its applicability. It also allows proofs to be obtained for an interval of parameters rather than for single parameter values.

Working with invariant manifolds and their intersections in our example would run into the following difficulties. If we considered time-shift maps, the system would need to be integrated in interval arithmetic over the whole parameter interval, and such integration quickly accumulates large errors; estimates on the manifolds would then have to be either local or restricted to very narrow parameter ranges. If instead we worked with the vector field in the extended phase space, we would need to establish the existence of a hyperbolic periodic orbit together with its stable and unstable manifolds, and handling a parameter interval would require parameterising these manifolds to include the parameter. These technicalities can be overcome, but only at the cost of difficulties to which our method is entirely immune. We emphasize that our method concludes the existence of a meaningful rotational
horseshoe, which is not any horseshoe which might just involve local chaos.

\medskip
\noindent\textbf{AI disclosure.}
Artificial intelligence was used solely to improve the
style and grammar of the manuscript and, after the
completion of the paper and its code, to perform an
independent audit of the implementation.
All mathematical ideas and arguments are entirely
the authors' own, and all code and its documentation
were written without AI assistance.
The subsequent audit, conducted using ChatGPT Astra,
found no errors or gaps and confirmed that the
implementation faithfully follows the methods
described in the paper.


\section{Prelimiaries}\label{s.preliminaries}

We now provide the details of the elements appearing in the statements introduced above.

\smallskip

All the theoretical content required for this work is topological in nature. 
Therefore, we may choose a model for the annulus \eqref{eq:annulus}.
We denote by $\boldsymbol{\pi}:\R^2\to \A$ the quotient map and by $\textrm{Homeo}_{\,0}(\A)$ the set of homeomorphisms in the homotopy class of the identity.

To define a \emph{rotational horseshoe} we consider a topological rectangle 
$\tilde{R} \subset \R^2$, whose boundary is given by the concatenation of arcs (see Figure \ref{f.horseshoe})

\[
\tilde{v}_1 \cdot \tilde{h}_1 \cdot \tilde{v}_2 \cdot \tilde{h}_2,
\]
and take $R:=\boldsymbol{\pi}(\tilde R)\subset \mathbb{A}$. We define integer horizontal translations $\tilde R_j \subset \mathbb{R}^2$ of $\tilde R$ as
\[\tilde R_j := \tilde R + (j 2\pi,0 )\qquad \mbox{for }j\in \mathbb{Z}.\]
(Note that $\boldsymbol{\pi} (\tilde R_j) = R$, for every $j\in \mathbb{Z}$.) We will say that $R$ is a \emph{rotational horseshoe} under the map $f$ if  for any lift $F$ of $f$ to the plane, $F(\tilde R)$ has Markovian intersections with at least two distinct integer horizontal translations $\tilde R_{j_0}$ and $\tilde R_{j_1}$, $j_0 \ne j_1$, in the following sense.

By a Markovian intersection between $F(\tilde R)$ and $\tilde R_j$, we mean a connected component of $F(\tilde R) \cap \tilde R_j$ whose boundary contains a sub-arc $\tilde h'_1$ of $F(\tilde h_1)$ and another sub-arc $\tilde h'_2$ of 
$F(\tilde h_2)$, with $\tilde h'_1$ and $\tilde h'_2$ connecting points of 
$v_1 + (j 2\pi,0)$ and $v_2 + (j 2\pi,0)$ inside $\tilde R_j$ and so that $F(\tilde v_i),i=1,2$ does not meet such connected component
; see Figure \ref{f.horseshoe} where we have Markovian intersections between $F(\tilde R)$ and $\tilde R_j$ for $j=-1,0,1$.

A configuration where we have a Markovian intersection between $F(\tilde R)$ and $\tilde R_j$ for $j\in J$, where $J$ is a finite subset of $\mathbb{Z}$, ensures a 
semiconjugacy on an invariant set $\Lambda \subset R$, which supports a map 
$h:\Lambda \to J^{\mathbb{Z}}$,
constructed exactly as Smale introduced in the hyperbolic case. 

Moreover, due to the separation properties of the plane, 
every periodic word in $J^{\mathbb{Z}}$  corresponds to a periodic point in $\Lambda$ (see \cite{GuiheneufMiliton2025}), and the word encodes the rotational information of that periodic point. For instance, a periodic point corresponding to the $n$-periodic word $j_1 j_2 \ldots j_n$, where $j_1,\ldots,j_n$ are (not necessarily distinct) elements of $J$, implies a translation from $\tilde R_0$ to $\tilde R_{j_1+\ldots+j_n}$ after a single period of the corresponding orbit. Moreover, the $j_0,\ldots,j_n$ encode the successive translations after each iterate of the map during the evolution of the periodic orbit. 
Finally, as a consequence of semicontinuity, the topological entropy of $\Lambda$ is bounded below by 
$\ln( |J|)$, where $|J|$ denotes the cardinality of $J$.


Associated with a rotational horseshoe is the \emph{instability region}. It can be defined as the smalles open annulus such that
it contains the unstable and stable sets associated to the horseshoe (constructed following the classical Smale construction).
It is easy to prove that such set exists, that it is an essential open annulus in $\A$ and that 
contains no a separating invariant essential continumm: that is a continuum
which separate the two ends of the annulus.

To any invariant set $\mathcal{K}$ such as the instability region, we can associate the so 
called \emph{rotation set} for a given lift $F$ of $f$.
Denoting by $\textrm{pr}_1$ the projection onto the first coordinate of $\R^2$, this set is given by

\begin{align*}\rho_\mathcal{K}(F) = \Bigl\{y\in \R \, \Bigl | \, &y=\lim_{k\to +\infty}\frac{\textrm{pr}_1(F^{n_k}(\tilde x_k)-\tilde x_k)}{n_k} \\ 
&\mbox{for some sequence } \tilde x_k\in \boldsymbol{\pi}^{-1}(\mathcal{K}) \mbox{ and } n_k\to+\infty 
  \Bigr \}.
  \end{align*}

For $f$ homotopic to the identity, which will be our setting in this paper, such set is given by an (possibly unbounded) interval. We say that a rational 
number $\frac{p}{q}$ 
from a rotation set $\rho_{\mathcal{K}}F$ is realized 
by some periodic point $x\in\mathcal{K}$ 
when
\[F^q(\tilde{x})=\tilde{x}+(p 2\pi,0),\]
for
every $\tilde x$ from $\pi^{-1}(x)$.


\section{Underlying Theoretical Results}\label{s.theoreticalres}

\label{s.theoretical}

The theoretical framework for the proof of Theorems~A and B is developed in
\cite{passpir}. Before recalling some of the main results contained therein,
we introduce the required notions, the first of which is crucial for the
detection of rotational horseshoes.

Given a pair of topological disks $U_{0}, U_{1}$, we want to introduce
its \emph{rotational difference} for a map $f \in\mathrm{Homeo}_{0}(\A)$, a
notion that is quite intuitive. To do so, we start with a definition of a
disjoint pair of disks.

\begin{definition}
We will say that $(U_{0},U_{1})$ is a \emph{disjoint pair of disks} for $f$ if
the following conditions are satisfied:

\begin{itemize}
\item $U_{0} \cap f(U_{0}) \neq\emptyset$ and $U_{1} \cap f(U_{1})
\neq\emptyset$;

\item the sets $U_{0} \cup f(U_{0})$ and $U_{1} \cup f(U_{1})$ are disjoint
and inessential.
\end{itemize}
\end{definition}

\begin{definition}
For a disjoint pair of disks $(U_{0},U_{1})$ we define their \emph{rotational
difference} with respect to $f$ as
\[
\rho_{f}(U_{0},U_{1}) = k_{1} - k_{0},
\]
where the integers $k_{0},k_{1}$ are determined by the identities
\[
F(\widetilde U_{i})\cap(\widetilde U_{i} + (k_{i} 2\pi,0 ))\neq\emptyset,
\qquad\mbox{for }i=0,1,
\]
for any lift $F$ of $f$ and any lifts $\widetilde U_{0}$ and $\widetilde
U_{1}$ of $U_{0}$ and $U_{1}$, respectively.
\end{definition}

It is easy to see that $\rho_{f}(U_{0},U_{1})$ does not depend on the choice
of $F,\widetilde U_{0}$ and $\widetilde U_{1}$. The last two notions we need
in order to state our main result are the following:

\begin{definition}
Given two open sets $U_{0},U_{1}\subset\A$, we say that \emph{$U_{0}$ visits
$U_{1}$} if some forward orbit of a point in $U_{0}$ meets $U_{1}$.
\end{definition}

\begin{definition}
We say that the ends of $\A$ are \emph{Birkhoff-related} for $f$, whenever for
any pair of neighborhoods $U_{+}$ of the end $+\infty$ of $\A$ and $U_{-}$ of
the end $-\infty$ of $\A$, it holds that $U_{-}$ visits $U_{+}$ and $U_{+}$
visits $U_{-}$.
\end{definition}

The following result will be the main tool for the proof of Theorem A. By
$\mathrm{Homeo}_{0,\lambda}(\A)$ we mean the space homeomorphisms preserving
an area form $\lambda$ in $\A$.

\begin{thm}[Theorem B of \cite{passpir}]\label{t.fromotherarticle}
Let $f \in \mathrm{Homeo}_{0,\lambda}(\A)$, and let
$U_0, U_1$ be a disjoint pair of discs for which
$$U_0\mbox{ visits }U_1\ ,\ U_1\mbox{ visits }U_0\ ,\ \mbox{ and }\rho_f(U_0,U_1) = \rho\geq 3.$$
Further, assume that either:
\begin{enumerate}
\item[(1)] The ends of $\A$ are Birkhoff-related,
\end{enumerate}
or
\begin{enumerate}
\item[(2)] There exist $V^{+}$ and $V^{-}$, neighborhoods of $+\infty$ and $-\infty$ respectively,
which do not intersect $U_0 \cup U_1$, and such that $V^{+}$ does not visit $V^{-}$
and $V^{-}$ does not visit $V^{+}$.
\end{enumerate}
Then, some power of $f$ carries a rotational horseshoe whose instability region
$\mathcal{I}$ meets both $U_0$ and $U_1$.
Moreover, for some lift $F$ of $f$, the rotation set of $\mathcal{I}$ verifies
\[
\rho_{\mathcal{I}}(F) \supset \left[1, \rho - 1\right],
\]
and every rational number in
$\rho_{\mathcal{I}}(F)$
is realized within~$\mathcal{I}$.
\end{thm}

\subsection{Application to Second Order ODE}

Based on Theorem \ref{t.fromotherarticle} we can now state a theoretical
result that serves as the foundation of a computer--assisted methodology for
proofs of chaos for differential equations, such as the Driven Pendulum. The
key feature here is that we consider differential equations on $\A$ which are
$T$-periodic, and have a time reversing symmetry, which is defined as follows.

\begin{definition}
We say that a differential equation on $\A$ is \emph{$1$-skew
time--reversible}\footnote{here the $1$ indicates the first variable in
$\mathbb{R}^{2}$.}, when the following holds: If $(\tilde{x}%
(t),\tilde{y}(t))$ is a lift of a solution, then $(-\tilde{x}%
(-t),\tilde{y}(-t))$ is also a lift of a solution.
\end{definition}

\begin{lemma}
\label{lem:symmetry}The driven pendulum (\ref{eq-driven-pendulum}) is $1$-skew time--reversible.
\end{lemma}

\begin{proof}
Assume that $(\tilde{x}(t),\tilde{y}(t))$ is a lift of a solution. We need to
show that $\left(  \tilde{p}(t),\tilde{q}(t)\right)  =(-\tilde{x}%
(-t),\tilde{y}(-t))$ satisfies the ODE (\ref{eq-driven-pendulum}). We directly
compute%
\[
\tilde{p}^{\prime}(t)=\frac{d}{dt}\left(  -\tilde{x}(-t)\right)  =\tilde
{x}^{\prime}(-t)=\tilde{y}(-t)=\tilde{q}(t),
\]
and%
\begin{multline*}
\tilde{q}^{\prime}(t)=\frac{d}{dt}\left(  \tilde{y}(-t)\right)  =-\tilde
{y}^{\prime}(-t)=-\left(  -\frac{g}{\ell}\sin(\tilde{x}\left(  -t\right)
)+A\sin(\omega_{0}\left(  -t\right)  )\right)  \\
=-\frac{g}{\ell}\sin(\tilde{p}(t))+A\sin(\omega_{0}t),
\end{multline*}
which concludes our proof.
\end{proof}

The reversible symmetry is reflected on the time map along the solution of the ODE. To state the result we introduce the notation $S:\mathbb{A\rightarrow A}$ for a map defined follows. For $\left(
x,y\right)  =\boldsymbol{\pi}\left(  \tilde{x},\tilde{y}\right)  \in
\mathbb{A}$ we define
\[
S\left(  x,y\right)  :=\boldsymbol{\pi}\left(  -\tilde{x},\tilde{y}\right)  .
\]
For $\left(  \tilde{x},\tilde{y}\right)  \in\mathbb{R}^{2}$ we will write
$\tilde{S}\left(  \tilde{x},\tilde{y}\right)  :=\left(  -\tilde{x},\tilde
{y}\right)  $.

\begin{corollary}\label{cor:S-symmetry}
Let $f:\mathbb{A\rightarrow A}$ be the time $T=\frac{2\pi}{\omega_{0}}$ shift
along the trajectory along the flow of the driven pendulum
(\ref{eq-driven-pendulum}) starting at the initial time equal to zero. Then%
\begin{equation}
f\circ S=S\circ f^{-1}.\label{eq:S-symmetry-maps}%
\end{equation}

\end{corollary}

\begin{proof}
For $\left(  x_{0},y_{0}\right)  \in\mathbb{A}$ let $\left(  x\left(
t\right)  ,y\left(  t\right)  \right)  $ be the solution starting from
$\left(  x_{0},y_{0}\right)  $ at the time $t=0$, i.e. $\left(  x\left(
0\right)  ,y\left(  0\right)  \right)  =\left(  x_{0},y_{0}\right)  $. Let
$\left(  \tilde{x}\left(  t\right)  ,\tilde{y}\left(  t\right)  \right)  $ be
a lift of the solution. Let us then write $\Phi_{t}\left(  x_{0},y_{0}\right)
:=\left(  \tilde{x}\left(  t\right)  ,\tilde{y}\left(  t\right)  \right)  $. 

From Lemma \ref{lem:symmetry} we know that%
\[
\Phi_{t}\left( \tilde S\left(  \tilde{x}_{0},\tilde{y}_{0}\right)  \right)
=\Phi_{t}(-\tilde{x}_{0},\tilde{y}_{0})=(-\tilde{x}\left(  -t\right)
,\tilde{y}\left(  -t\right)  )= \tilde S(\tilde{x}\left(  -t\right)  ,\tilde
{y}\left(  -t\right)  )= \tilde S\left(  \Phi_{-t}\left(  \tilde{x}_{0},\tilde
{y}_{0}\right)  \right)  .
\]
Taking $t=T$ and writing $\tilde f$ for $\Phi_{T}$, we obtain
\[
\tilde{f}\circ\tilde{S}=\Phi_{T}\circ\tilde{S}=\tilde{S}\circ\Phi_{-T}%
=\tilde{S}\circ\tilde{f}^{-1},
\]
which implies (\ref{eq:S-symmetry-maps}).
\end{proof}

The proof of the following theorem relies on Theorem \ref{t.fromotherarticle}
and standard arguments from the field of surface dynamics.

\begin{thm}\label{t.diffeq}
Consider a $1$-skew time--reversible, periodically forced differential equation in $\A$, whose associated Poincar\'e map
$f$ preserves the Lebesgue measure.
Assume that for some disjoint pair of disks $(U_0,U_1)$ for $f$
it holds that
$$\ U_0\mbox{ visits }U_1, \, U_1\mbox{ visits }U_0 \mbox{ and } \rho_f(U_0,U_1)\geq 3.$$
Then, some power of $f$ has a rotational horseshoe.
Moreover, the associated instability region $\mathcal{I}$
meets $U_0,U_1$ and
the rotation set for any lift $F$ of $f$ 
contains an interval of length $\rho_f(U_0,U_1)-2$. In particular, every rational point inside such interval is realized by some periodic orbit inside
$\mathcal{I}$.
\end{thm}

\begin{proof}
Consider the following two cases.

\begin{enumerate}
\item The ends of $\A$ are Birkhoff-related.

\item The ends of $\A$ are not Birkhoff-related.
\end{enumerate}

In the first case the result follows directly from Theorem
\ref{t.fromotherarticle}. What remains is to treat the second case.

Since the ends are not Birkhoff-related, without the loss of generality we can
assume that for some $r>0$ and for
\[
V^{+}=\pi\left(  \{(\tilde x, \tilde y):\, \tilde y\geq r\}\right)
\mbox{\quad and \quad}V^{-}=\pi\left(  \{(\tilde x, \tilde y):\, \tilde
y\leq-r\}\right)
\]
no point in $V^{-}$ visits $V^{+}$ in the future under iterations of the map
$f$. (We also have the second possibility where no point in $V^{+}$ visits
$V^{-}$, which can be treated by mirror arguments). We can choose $r>0$ as
large as we want, so we take $r$ large enough so that $V^{-}$ and $V^{+}$ are
disjoint from $U_{0}$ and $U_{1}$.

If for some time $t>0$ we had a lift of some solution starting from $(\tilde
x(0),\tilde y(0)) \in\pi^{-1}(V^{+})$, which would reach $(\tilde x(t),\tilde
y(t)) \in\pi^{-1}(V^{-})$, then, by the $1$-skew time--reversibility of the
system, the solution starting from the point $(-\tilde x(t),\tilde y(t))\in
\pi^{-1}(V^{-})$ would reach $(-\tilde x(0),\tilde y(0))\in\pi^{-1}(V^{+})$
after the time $t>0$. This would contradict the fact that $V^{-}$ does not
visit $V^{+}$ in the future; hence $V^{+}$ cannot visit $V^{-}$ in the future.
We are therefore in the second condition of Theorem \ref{t.fromotherarticle},
hence we are done.
\end{proof}

\section{Strategy}

\label{sec-strategy}

In this section we fix parameters $g$, $\ell$, $A$, and $\omega_{0}$ in
\eqref{eq-driven-pendulum}, and let $F:\mathbb{R}^{2}\rightarrow\mathbb{R}%
^{2}$ denote the induced time-$T$ Poincaré \ map, where
$T=2\pi/\omega_{0}$. In other words, for any initial condition $x\in
\mathbb{R}^{2}$, $F(x)$ is the state of the corresponding solution at time
$T$. If $f:\mathbb{A}\rightarrow\mathbb{A}$ is the associated Poincar\'{e} map
on the annulus, then $F$ is a lift of $f$.

Our goal is to determine whether, for the given parameters, the hypotheses of
\textit{Theorem~\ref{t.diffeq}} are satisfied. To this end, we numerically
look for a pair of sets $\tilde{R}_{0}$ and $\tilde{R}_{1}$ homeomorphic to
two discs in $\mathbb{R}^{2}$, which will serve as our \emph{candidates}, and
rigorously verify that the conditions given in \textit{Corollary~\ref{candidates}%
} are satisfied. The corollary follows from \textit{Theorem~\ref{t.diffeq}}.

\begin{corollary}
\label{candidates} Let $\tilde{R}_{0}$ and $\tilde{R}_{1}$ be homeomorphic to two discs in
$\mathbb{R}^{2}$, such that $U_{0}:=\boldsymbol{\pi} \tilde{R}_{0}$ and $U_{1}:=\boldsymbol{\pi} \tilde
{R}_{1}$ are homeomorphic to two discs in $\A$. Suppose that, for $i=0,1$, there exist
$k_{i},l_{i}\in\mathbb{Z}$ and $N_{i}\in\mathbb{N}$, $N_i\ge 2$, such that the
following conditions hold: 

\begin{enumerate}

\item $\left(  \tilde{R}_{i}+(k_{i}2\pi,0)\right)  \cap F(\tilde{R}_{i}%
)\neq\emptyset$. 

\item $\left(  \tilde{R}_{i}+((k_{i}\pm1)2\pi,0)\right)  \cap F(\tilde{R}%
_{i})=\emptyset$. 

\item $U_{i}\cap U_{1-i}=\emptyset$ and $f(U_{i})\cap U_{1-i}=\emptyset$.

\item $|k_{1} - k_{0}| \ge3$. 

\item $f^{N_{i}}(U_{i}) \cap U_{1-i} \neq\emptyset$. 
\end{enumerate}

Then $U_{0}$ and $U_{1}$ satisfy the hypotheses of {Theorem~\ref{t.diffeq}}, hence the claims of the theorem hold true.
\end{corollary}

In {Section~\ref{section-implementation}} we explain how the above conditions
can be rigorously verified using CAPD, even though no explicit formula for $F$
is available.

Below we describe the choice of sets used in the proof of
Theorem~A.

The sets $U_{0}$ and $U_{1}$ can be found by searching for four points
$x_{0},\ldots,x_{3}\in\mathbb{A}$ such that $x_{0},x_{2}\in U_{0}$,
$x_{1},x_{3}\in U_{1}$ and
\[
f\left(  x_{0}\right)  \in U_{0},\qquad f\left(  x_{1}\right)  \in
U_{1},\qquad f^{N_{0}}\left(  x_{3}\right)  \in U_{0},\qquad f^{N_{1}}\left(
x_{2}\right)  \in U_{1}.
\]

Intuitively, we can search for points such that%
\begin{align*}
&  f(x_{0})\text{ is close to }x_{2},\\
&  f\left(  x_{1}\right)  \text{ is close to }x_{3},\\
&  f^{N_{0}}\left(  x_{3}\right)  \text{ is close to }x_{0},\\
&  f^{N_{1}}\left(  x_{2}\right)  \text{ is close to }x_{1}.
\end{align*}
In practice, we can search for $\tilde{x}_{0},\ldots,\tilde{x}_{3}%
\in\mathbb{R}^{2}$, by minimising%
\begin{align}
&  \left\Vert F\left(  \tilde{x}_{0}\right)  -\left(  \tilde{x}_{2}+\left(
2k_{0}\pi,0\right)  \right)  \right\Vert ^{2}+\left\Vert F\left(  \tilde
{x}_{1}\right)  -\left(  \tilde{x}_{3}+\left(  2k_{1}\pi,0\right)  \right)
\right\Vert ^{2}\label{eq:x0-x3-choice}\\
&  +\left\Vert F^{N_{0}}\left(  \tilde{x}_{3}\right)  -\left(  \tilde{x}%
_{0}+\left(  2l_{0}\pi,0\right)  \right)  \right\Vert ^{2}+\left\Vert
F^{N_{1}}\left(  \tilde{x}_{2}\right)  -\left(  \tilde{x}_{1}+\left(
2l_{1}\pi,0\right)  \right)  \right\Vert ^{2},\nonumber
\end{align}
where $k_{0},k_{1},l_{0},l_{1}\in\mathbb{Z}$ are some chosen numbers; then
$x_{i}=\boldsymbol{\pi} \tilde{x}_{i}$ for $i=0,\ldots,3$.

\begin{remark}
In our application we will take $k_{0},k_{1}$ so that $\left\vert k_{0}%
-k_{1}\right\vert \geq3$ to ensure condition (4) from Corollary \ref{candidates}.
\end{remark}

\begin{figure}[ptb]
\begin{center}
\includegraphics[width=8cm]{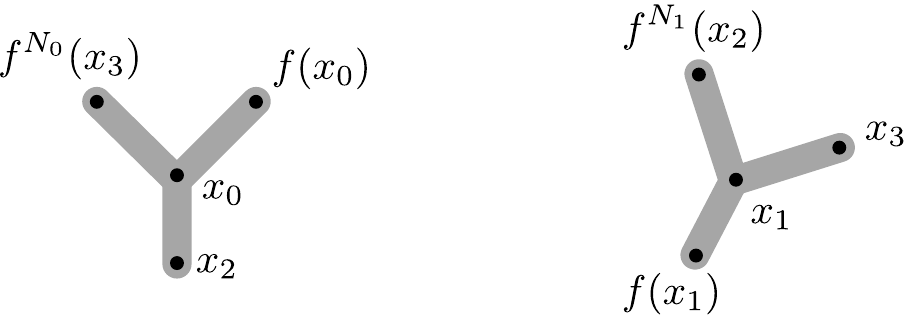}
\end{center}
\caption{Two branched-sets $U_{0}$ (left) and $U_{1}$ (right), in grey. Each
of the two sets has three branches.}%
\label{fig:y-sets}%
\end{figure}

Once $x_{0},\ldots,x_{3}\in\mathbb{A}$ are found, we can take $U_{0}$ and
$U_{1}$ as sets which contain these points, as depicted in Figure
\ref{fig:y-sets}. We refer to such sets as \textquotedblleft branched-sets".

\begin{remark}
\label{rem:cond-5}With such choice of sets, since $x_{3}\in U_{1}$, $f^{N_{0}%
}(x_{3})\in U_{0}$ and $x_{2}\in U_{0}$, $f^{N_{1}}(x_{2})\in U_{1}$, we are
sure that condition (5) from Corollary \ref{candidates} will be fulfilled.
\end{remark}

\begin{figure}[ptb]
\begin{center}
\includegraphics[width=8cm]{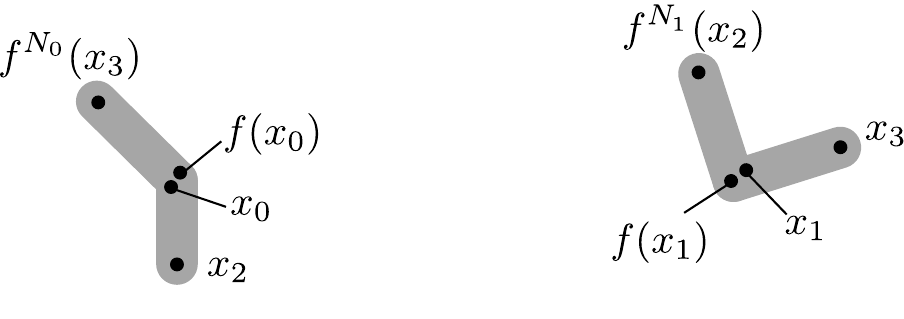}
\end{center}
\caption{Two branched-sets $U_{0}$ (left) and $U_{1}$ (right), each with two
branches.}%
\label{fig:v-sets}%
\end{figure}

When $x_{0}$ and $x_{1}$ are close to two fixed points, meaning that
$x_{0}\approx f\left(  x_{0}\right)  $ and $x_{1}\approx f\left(
x_{1}\right)  $, then we can take branched-sets of the shape from Figure
\ref{fig:v-sets}.



\begin{figure}[ptb]
\begin{center}
\includegraphics[width=7cm]{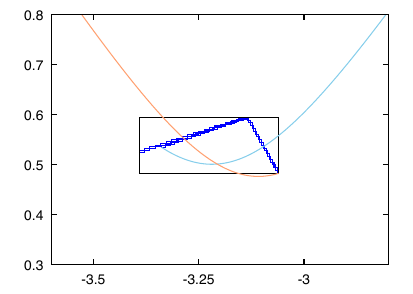}\includegraphics[width=7cm]{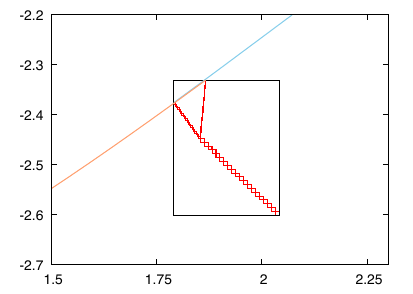}
\end{center}
\caption{The choice of two branched-sets $\tilde R_{0}$ (in blue) and $\tilde R_{1}$ (in
red). }%
\label{fig:Ri-sets}%
\end{figure}

\begin{figure}[ptb]
\begin{center}
\includegraphics[width=14cm]{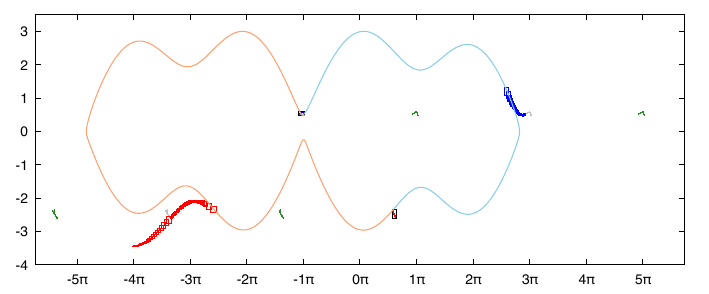}
\end{center}
\caption{The bounds on the image of the branched-sets $\tilde R_{0}$ and $\tilde R_{1}$ in
blue and red, respectively. The light blue curve depicts the trajectory along
the flow of the ODE from $\tilde R_{0}$ to $\tilde R_{1}$. In light red we depict the
trajectory from $\tilde R_{1}$ to $\tilde R_{0}$. }%
\label{fig:FRi-sets}%
\end{figure}

\begin{figure}[ptb]
\begin{center}
\includegraphics[width=14cm]{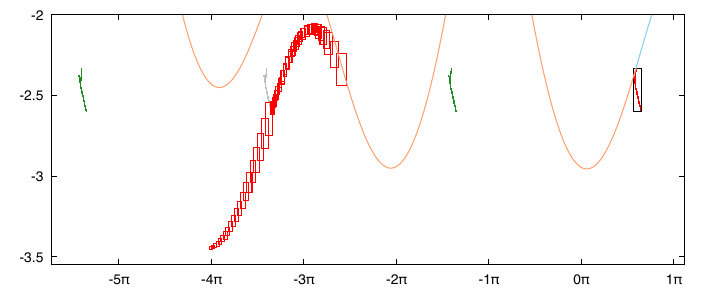}
\end{center}
\caption{The bounds on the image of the branched-set $\tilde R_{1}$ in red. In grey
we depict $\tilde R_{1}+(k_{1} 2\pi,0)$. (We take $k_{1}=-2$.) In green we depict
$\tilde R_{1}+((k_{1}\pm1) 2\pi,0)$. }%
\label{fig:FR1-set}%
\end{figure}

The branched-sets which we consider for the proof of Theorem~A are depicted in Figure \ref{fig:Ri-sets}. In Figure
\ref{fig:FRi-sets} we depict the bounds on their images, together with the
trajectories along the flow of the ODE. These trajectories demonstrate that
the sets visit each other, since the trajectories start at $\tilde{R}_{i}$ and
end in $\tilde{R}_{1-i}$, for $i=0,1$. Figure \ref{fig:FR1-set} gives a
closeup of the image of the set $\tilde{R}_{1}$ together with $\tilde{R}%
_{1}+(k_{1}2\pi,0)$ and $\tilde{R}_{1}+((k_{1}\pm1)2\pi,0)$.

\subsection{CAPD implementation}

\label{section-implementation}

The verification of the conditions in {Corollary~\ref{candidates}} is carried out
using interval arithmetic and computer-assisted proofs. Rather than working
with single floating-point values, all computations are performed using
intervals, which represent sets of real numbers and are propagated in such a
way that all possible rounding and truncation errors are rigorously accounted
for. As a result, the outcome of each operation encloses the set of all
possible exact results.

In our implementation, interval arithmetic is provided by the CAPD library.
Within this framework, the basic data type \texttt{interval} represents a
closed interval in $\R$. When working in dimension two, this notion is
naturally extended to the data type \texttt{IVector}, which represents a
vector of intervals and can be interpreted as a rectangular subset of $\R^{2}$.

The CAPD library also provides tools for the rigorous integration of ordinary
differential equations with interval-valued initial conditions. In particular,
one can define an \texttt{ITimeMap}, which represents the time-$T$ flow of a
vector field whose coefficients are treated using interval arithmetic. When an
\texttt{ITimeMap} is applied, the result is another
\texttt{IVector} that is guaranteed to contain the exact image of the initial
set under the flow at the specified time.

In the present work, this framework is used to rigorously enclose the image of
points or rectangular subsets of $\R^{2}$ under the flow of the driven
pendulum equation.

When working with CAPD we perform all our computations by integrating the ODE
in $\mathbb{R}^{2}$. A box in $\mathbb{A}$ which corresponds to a given box in
$\mathbb{R}^{2}$ is obtained by a shift modulo $2\pi$ along its first
coordinate. We can easily check that two boxes in $\mathbb{R}^{2}$ are
disjoint. We can also check that their projections onto $\mathbb{A}$ are
disjoint by considering their modulo $2\pi$ shifts.


\begin{proof}
[Proof of Theorem A]We consider $F:\R^{2}\rightarrow
\R^{2}$ as the time-$T$ Poincar\'{e} map induced by
\eqref{eq-driven-pendulum}, with $T=2\pi/\omega_{0}$.

We have chosen $k_{0}=2,k_{1}=-2,N_{0}=3,N_{1}=2$ and $l_{0}=l_{1}=0$. For
this choice we have found points $\tilde{x}_{0},\ldots,\tilde{x}_{3}%
\in\mathbb{R}^{2}$ which minimise (\ref{eq:x0-x3-choice}) by non-rigorous
numerical investigation. This resulted in
\begin{align*}
\tilde{x}_{0} &  =(3.1487248116739,0.59522282432838),\\
\tilde{x}_{1} &  =(1.8534555698055,-2.4470505318419),\\
\tilde{x}_{2} &  =(2.9447006778146,0.53538184902088),\\
\tilde{x}_{3} &  =(1.866868002463,-2.3302174559969).
\end{align*}
We have then computed rigorous-interval-arithmetic bounds on the images of
$F(\tilde{x}_{0})$, $F(\tilde{x}_{1})$, $F^{N_{0}}(\tilde{x}_{3})$ and
$F^{N_{1}}(\tilde{x}_{2})$ and shifted these sets modulo $2\pi$ along the
first coordinate to choose our branched-sets $\tilde{R}_{0},\tilde{R}%
_{1}\subset\mathbb{R}^{2}$; see Figure \ref{fig:Ri-sets}. We have made sure
that we shift the sets which enclose $F(\tilde{x}_{i})$ by $k_{i}2\pi$, which
ensures condition (1) from Corollary \ref{candidates}.

The two black rectangles in Figure \ref{fig:Ri-sets} depict the two boxes
(\ref{eq:th-boxes}) from the statement of Theorem~A. (In
the Figures \ref{fig:Ri-sets}, \ref{fig:FRi-sets} and \ref{fig:FR1-set} we
have shifted the box enclosing $\tilde{R}_{0}$ to the left by $2\pi$. This was
done to make it clearly visible on the plots that the two sets visit each other.)

By Remark \ref{rem:cond-5}, our branched-sets ensure that condition (5) from
Corollary \ref{candidates} is fulfilled.

We have computed the bounds on the images $f(U_{i})$ by subdividing each
branch of the sets into $20$ subsets\footnote{Subdividing each branch into
$10$ subsets would suffice, but this would result in less transparent plots
than those in Figures \ref{fig:Ri-sets}, \ref{fig:FRi-sets} and
\ref{fig:FR1-set}.}, and integrating the ODE in interval arithmetic. We have
used these bounds to validate conditions (2) and (3) from Corollary
\ref{candidates}.

To provide an illustration, the bounds from our computer assisted proof are
included in Figures \ref{fig:Ri-sets}, \ref{fig:FRi-sets} and
\ref{fig:FR1-set}. The plots clearly indicate that all conditions from Corollary
\ref{candidates} are fulfilled.

The code for the proof is available at \cite{github} and the computation time
is under one second on a standard laptop.
\end{proof}

\section{Parameters Exploration}

In this section we show how the process of choosing branched sets can be
automated. To do so we exploit the fact that the driven pendulum has periodic
orbits, close to which we can position our branched sets. We will also make
use of the $S$-symmetry of the system.

Branched sets are constructed based on the choice of four points $x_{0},\ldots,x_{3}$, as depicted in Figures \ref{fig:y-sets} and
\ref{fig:v-sets}. One way to find these points is by minimising
(\ref{eq:x0-x3-choice}). This does not always work well in practice though. If the
system has periodic orbits then we can choose branched sets of the form 
Figure \ref{fig:v-sets}, and instead of minimising (\ref{eq:x0-x3-choice}) 
perform the following two steps:

\begin{enumerate}
\item Find $x_{0}$ and $x_{1}$, which are close to some
periodic orbits.

\item Find $x_{2}$ close to $x_{0},$ which after $N_{1}$
iterates lands close to $x_{1}$; find $x_{3}$ close to
$x_{1}$, which after $N_{0}$ iterates lands close to $x_{0}$.
\end{enumerate}

The first step is best performed in the lift by choosing some $k_0,k_1\in \mathbb{Z}$, $|k_0-k_1|\ge3$, and solving%
\begin{equation}
F(\tilde{x}_{i})-(\tilde{x}_{i}+(k_{i}2\pi,0))=0\qquad\text{for }%
i=0,1.\label{eq:find-zeros}%
\end{equation}
It is not necessary to validate that$\ \tilde{x}_{0}$ and $\tilde{x}_{1}$ are
true solutions of (\ref{eq:find-zeros}). It is enough that they are 
numerical approximations. We then take $x_{i}=\boldsymbol{\pi}\tilde{x}_{i}%
\in\mathbb{A}$, for $i=0,1$ . 

For the second step one can consider many randomly sampled points $\hat{x}%
\in\mathbb{A}$ close to $x_{0}$, and check if $f^{n}\left(  \hat{x}\right)  $,
for some $n\geq2$, is close to $x_{1}$. We then choose $N_{1}:=n$ and
$x_{2}:=\hat{x}$ for the point $\hat{x}$ for which $f^{n}\left(  \hat
{x}\right)  $ lands closest to $x_{1}$. The point $x_{3}$ can be chosen in a
similar way. 

\begin{figure}
\begin{center}
	\includegraphics[width=10cm]{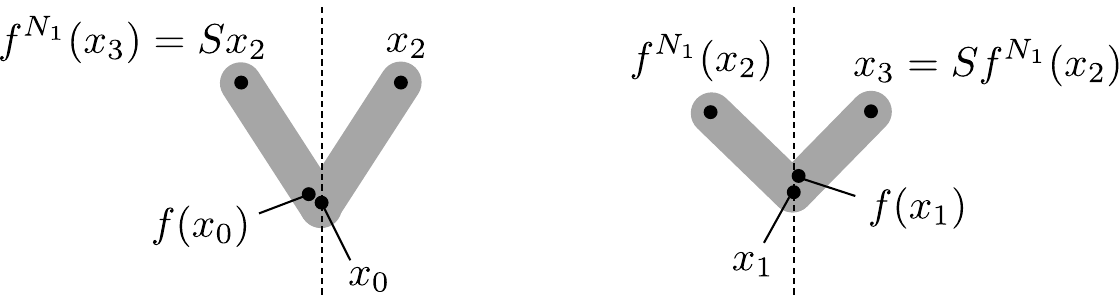}
\end{center}
\caption{Two branched sets $U_0$ on the left and $U_1$ on the right, each with with two branches chosen from the $S$-symmetry. The dotted lines depict the set of self $S$-symmetric points. \label{fig:Sv-sets-fig} }
\end{figure}

In our application, instead of choosing $x_3$ in a similar way as $x_2$, we can make use of the fact that the
system is $S$-symmetric (see Corollary \ref{cor:S-symmetry}). It turns out that in step (1) we can choose the points $x_{0}$ and $x_{1}$
to be self $S$-symmetric. In this case, once $x_2$ is found, a good choice of $x_{3}$ is
$x_{3}:=Sf^{N_{1}}\left(  x_{2}\right)  $ and a good choice of $N_0$ is $N_{0}:=N_{1}$. Below lemma
justifies why this is a desirable choice (see also Figure \ref{fig:Sv-sets-fig}).

\begin{lemma}
If $x_{3}=Sf^{N_{1}}\left(  x_{2}\right)  $ then $Sx_{2}=f^{N_{1}}\left(
x_{3}\right)  $.
\end{lemma}

\begin{proof}
From the $S$-symmetry of $f$ (i.e. from $f^{-1}\circ S=S\circ f$; see Corollary \ref{cor:S-symmetry}), for every
$k\in\mathbb{Z}$ we have%
\[
f^{k}\circ S\circ f^{k}=f^{k}\circ(S\circ f)\circ f^{k-1}=f^{k}\circ
(f^{-1}\circ S)\circ f^{k-1}=f^{k-1}\circ S\circ f^{k-1}.
\]
Iterating this argument leads to%
\[
f^{N_{1}}\left(  x_{3}\right)  =f^{N_{1}}\left(  Sf^{N_{1}}\left(
x_{2}\right)  \right)  =f^{N_{1}}\circ S\circ f^{N_{1}}\left(  x_{2}\right)
=...=S x_{2},
\]
as required.
\end{proof}

Once the points $x_{0},x_1,x_2$ and $x_{3}=Sf^{N_1}(x_2)$ are selected as described
above, we choose two sets $U_{0}$ and $U_{1}$, with  $S$-symmetric branch
points, as depicted in Figure \ref{fig:Sv-sets-fig}.

The above described two steps (1) and (2) are easy to automate. With this procedure we have proven Theorem B.
\begin{proof}[Proof of Theorem B] We have subdivided the interval $[1,1.5]$ into $25\cdot 10^3$ sub-intervals and on each of them we have constructed branched sets $U_1$ and $U_2$ as described in the discussion preceding the proof. 

To solve (\ref{eq:find-zeros}) we have always chosen $k_0=-1$ and $k_1=2$. This is an arbitrary choice made by us. By trial and error we have found  that such choice produces the desirable branched sets in most cases.

When searching for $x_2$ we have allowed our automated numerical search to look for $N_1\le 4$. This was an important choice since for $n>4$ we have found that the interval enclosure for $f^{n}(\hat x)$ starts to blow up. This is because the interval integrator needs to work with an interval of parameters $A$. We take our $A$-intervals with a tiny $2\cdot 10^{-5}$ width, but even this produces a blowup when integrating too long. Our restriction that $N_1\le 4$ results in the longest allowed ODE integration time of $4\cdot 2\pi$.

Once $U_0$ and $U_1$ are chosen, we validate the conditions from Corollary \ref{candidates}. To compute the bounds for $f(U_0)$ and $f(U_1)$ we subdivide each branch of the sets into $20$ fragments. Again, such choice is quite arbitrary, but works in most cases.

\begin{figure}[ptb]
\begin{center}
\includegraphics[width=14cm]{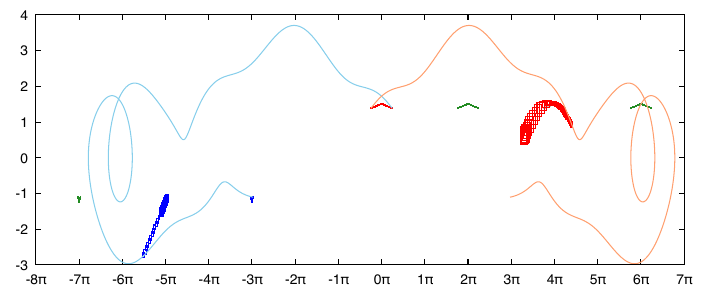}
\end{center}
\caption{The bounds on the image of the branched-sets $\tilde R_{0}$ and $\tilde R_{1}$ in
blue and red, respectively. The light blue curve depicts the trajectory along
the flow of the ODE from $U_{0}$ to $ U_{1}$. In light red we depict the
trajectory from $U_{1}$ to $ U_{0}$. (The light red trajectory does not lead from $\tilde R_1$ to $\tilde R_0$ in the lift, but it does lead to $\tilde R_0 + (6\pi,0)$, which is the same set on the annulus.) }%
\label{fig:FRi-sets-V}%
\end{figure}

In order to provide some intuition, in Figure \ref{fig:FRi-sets-V} we give an illustration of what happens for $A\in [1.25,1.25002]$. We have chosen this parameter range since it is in the middle of our considered interval $[1,1.5]$. For other parameters the pictures would be similar.

The computer assisted proof took under one and a half hours, running on a single thread on a standard laptop. The code is made available in \cite{github}. 
\end{proof}

\end{document}